\documentclass[11pt]{article}

\usepackage{mathrsfs}
\usepackage{latexsym}
\usepackage{amsmath,amsthm}
\usepackage{amsfonts}
\usepackage{url}

\usepackage{amssymb}
\newcommand{\h}{\eta}
\newcommand{\etab}{\bar\eta}

\usepackage{bbm}
\newcommand{\chr}{\boldsymbol{\mathbbm{1}}} % characteristic function
\newcommand{\pred}[1]{\chr_{\left\{ #1 \right\}}}
\newcommand{\TV}[1]{\nrm{#1}_{\textrm{{\tiny \textup{TV}}}}}

\newcommand{\Lip}[1]{\nrm{#1}_{\textrm{{\tiny \textup{Lip}}}}}

\newcommand{\xoi}{\set{0,1}}

\usepackage{geometry}
\makeatletter
\usepackage{mytheo}

\newcommand{\ben}{\begin{enumerate}}
\newcommand{\een}{\end{enumerate}}
\newcommand{\bit}{\begin{itemize}}
\newcommand{\eit}{\end{itemize}}

\newcommand{\basicspace}{\Omega}
\newcommand{\X}{\basicspace}
\newcommand{\supr}[1]{^{(#1)}}
\newcommand{\seq}[3]{(#1_{#2},\ldots,#1_{#3})}
\newcommand{\sseq}[3]{#1_{#2}^{#3}}  % short seq

\newcommand{\dsabs}[1]{\bigl| #1 \bigr|}
\newcommand{\nrm}[1]{\left\Vert #1 \right\Vert}

\newcommand{\calL}{\mathcal{L}}

\newcommand{\calP}{\mathcal{P}}

\newcommand{\calX}{\mathcal{X}}
\newcommand{\calY}{\mathcal{Y}}

\newcommand{\R}{\mathbb{R}}
\newcommand{\N}{\mathbb{N}}

\newcommand{\beq}{\begin{eqnarray*}}
\newcommand{\eeq}{\end{eqnarray*}}
\newcommand{\beqn}{\begin{eqnarray}}
\newcommand{\eeqn}{\end{eqnarray}}
\newcommand{\paren}[1]{\left( #1 \right)}
\newcommand{\sqprn}[1]{\left[ #1 \right]}
\newcommand{\tlprn}[1]{\left\{ #1 \right\}}
\newcommand{\set}[1]{\tlprn{#1}}
\newcommand{\abs}[1]{\left| #1 \right|}

\newcommand{\gn}{\, | \,}

\newcommand{\ds}{\displaystyle}
\newcommand{\ts}{\textstyle}

\renewcommand{\th}{\ensuremath{^{\mathrm{th}}}~}

\newcommand{\hide}[1]{}
\newcommand{\oo}[1]{\frac{1}{#1}}

\def\eps{\varepsilon}
\newcommand{\defeq}{=}
\newcommand{\serp}{\oplus}
\newcommand{\parp}{\otimes}

\title{Constructing processes with prescribed mixing coefficients}
\author{Leonid (Aryeh) Kontorovich\thanks{Supported in part by the Israel Science Foundation}\\
Department of Mathematics \\
Weizmann Institute of Science \\
Rehovot, Israel }

\begin{document}
\maketitle
\begin{abstract}
The rate at which dependencies between future and past observations decay in a random process may be quantified in terms of
mixing coefficients. The latter in turn appear in strong laws of large numbers and concentration of measure results for dependent
random variables. Questions regarding 
what rates are possible
%the possible mixing rates 
for various notions of mixing have been posed since the 1960's, and have
important implications for some open problems in the theory of strong mixing conditions.

This paper deals with $\eta$-mixing, a notion defined in [Kontorovich and Ramanan], which is closely related to $\phi$-mixing. We show that there
exist measures on finite sequences with essentially arbitrary $\eta$-mixing coefficients, as well as processes with arbitrarily slow
% or fast
mixing rates.
\end{abstract}

\section{Introduction}
\subsection{Preliminaries}

Strong mixing conditions deal with quantifying the 
%(usually decaying) 
decaying
dependence between blocks of random variables in 
a stochastic process. 
These have been traditionally used to establish strong laws of large numbers for non-independent processes.
Bradley \cite{bradley07i,bradley07ii,bradley07iii} is an encyclopedic source on the matter; see also his survey 
paper \cite{bradley05}. 
In \cite[
%p. 64
Chapter 26
]{bradley07iii},
Bradley  traces the early research on mixing rates to 
Volkonski{\u\i} and Rozanov \cite{volk61} and gives a comprehensive account of the progress since then.

Our interest in strong mixing was motivated by 
the desire for
concentration of measure bounds for non-independent random sequences. 
Given the excellent survey papers and monographs dealing with
concentration of measure (in particular, \cite{ledoux01},
\cite{lugosi03},
and
\cite{schechtman99}), we will give only the briefest summary here.

Suppose $\X$ is a finite\footnote{
The results hold verbatim for countable sets, and extend naturally to $\R$ under mild assumptions; see 
\cite{kont07-thesis,kont07-slln}.
} set and let $\mu$ be an arbitrary
(nonproduct) probability measure on $\X^n$. 
We proceed to define 
%the notion 
a type of
strong mixing used throughout this note.
%Let us define a type of mixing closely related to $\phi$-mixing.
For $1\leq i<j\leq n$
%, $y\in\X^{i-1}$ and 
%$w\in\X$, 
and $x\in\X^i$,
let
$$\calL(\sseq{X}{j}{n}\gn \sseq{X}{1}{i}=x)
$$ be the distribution
of $\sseq{X}{j}{n}\equiv\seq{X}{j}{n}$ conditioned on 
$\sseq{X}{1}{i}=x$. 
For $y\in\X^{i-1}$ and 
$w,w'\in\X$, 
define
\beqn
\label{eq:hdef}
\eta_{ij}(y,w,w') &=&
\TV{
\calL(\sseq{X}{j}{n}\gn \sseq{X}{1}{i}={y w})-
\calL(\sseq{X}{j}{n}\gn \sseq{X}{1}{i}={y w'})
},
\eeqn
where $\TV{\cdot}\equiv\oo2\nrm{\cdot}_1$ is the total variation norm;
%and
likewise,
define
\beqn
\label{eq:hbar}
\bar\eta_{ij} &=&
\max_{y\in\X^{i-1},w,w'\in\X}
\eta_{ij}(y,w,w').
\eeqn

%The 
%notion of mixing we define here 
This notion of mixing
is by no means new; it can be
traced (at least implicitly) to Marton's work \cite{marton98} and is
quite explicit in Samson \cite{samson00} and 
Chazottes et al. \cite{chazottes07}.
We are not aware of a
standardized term for this type of mixing, and have referred to it as
{\it $\eta$-mixing} in previous work \cite{kontram06}. It was observed in 
\cite{samson00} that the $\phi$-mixing coefficients bound the $\eta$-mixing ones:
\beq
\bar\h_{ij}&\leq& 2\phi_{j-i},
\eeq
and conjectured in \cite{kont07-thesis} that
\beq
\label{eq:fhconj}
\oo2\sum_{i=1}^{n-1} \phi_i
&\leq&
1+\max_{1\leq i<n}\sqprn{\sum_{j=i+1}^n \etab_{ij}};
\eeq
the latter remains open.

In all instances, $\eta$-mixing has come up in the context of concentration of measure.
In particular, define
%Define 
$\Gamma$ and $\Delta$ to be upper-triangular $n\times n$ matrices,
with $\Gamma_{ii}=\Delta_{ii}=1$ and
\beqn
\label{eq:GDdef}
\Gamma_{ij} = \sqrt{\etab_{ij}},
\qquad
\Delta_{ij} = \etab_{ij}
\eeqn
for $1\leq i<j\leq n$.

Samson \cite{samson00} proved that 
any distribution $\mu$ on $[0,1]^n$ and any convex
$f:[0,1]^n\to\R$ with $\Lip{f}\leq 1$ (with respect to $\ell_2$) satisfy
\beqn
\label{eq:samson}
\mu\set{\abs{f-\mu f}>t} &\leq& 2\exp\paren{-\frac{t^2}{2\nrm{\Gamma}_2^2}}
\eeqn
where $\nrm{\Gamma}_2$ is the $\ell_2$ operator norm.

Chazottes et al. \cite{chazottes07}
and independently,
the author with K. Ramanan \cite{kontram06} showed that 
any distribution $\mu$ on $\X^n$ and any
$f:\X^n\to\R$ with $\Lip{f}\leq n^{-1/2}$ (with respect to the 
Hamming metric) satisfy
\beqn
\label{eq:kontram}
\mu\set{\abs{f-\mu f}>t} &\leq& 2\exp\paren{-\frac{t^2}{2\nrm{\Delta}_\infty^2}}
\eeqn
where $\nrm{\Delta}_\infty$ is the $\ell_\infty$ operator norm
($\nrm{\Delta}_\infty$ may be replaced by $\nrm{\Delta}_2$ and
\cite{chazottes07} achieves
a better constant in the exponent).

Results of type (\ref{eq:samson}) and (\ref{eq:kontram}) are known as {\em concentration of measure} inequalities; broadly, they assert
that
any ``sufficiently continuous'' function is tightly
concentrated about its mean. 
Such bounds have a remarkable range of applications, spanning abstract fields such as asymptotic Banach space theory
\cite{ball97,schechtman99} as well as more practical ones such as randomized algorithms 
\cite{dubhashi98} and machine learning \cite{boucheron05}.
% and statistics 
Strong laws of large numbers are readily obtained from concentration bounds 
\cite{kont07-slln}. 
%well-behaved function is 
%essentially constant (i.e., very close to its mean) almost everywhere.
%Kontorovich's thesis 
%\cite{kont07-thesis}
%presents various concentration of measure results for $\eta$-mixing processes.

Having motivated the study of mixing and measure concentration,
%Having established that mixing and concentration are amply worthy of study, 
let us turn to the behavior of the
$\eta$-mixing coefficients. It is immediate from the construction that $\etab_{ij}$ is an upper-triangular
$n \times n$ matrix satisfying
\bit
\item[(P1)] $\etab_{ij}=0$ for $i\geq j$
\item[(P2)] $0\leq \etab_{ij}\leq  1$ for $1\leq i<j\leq n$.
\eit
It is also 
%easy 
simple
to show (as we shall do below in Lemma \ref{lem:hijdecr}) that
\bit
\item[(P3)] $\bar\eta_{ij_2}\leq\bar\h_{ij_1}$ for 
$
%1\leq 
i<j_1 < j_2
%\leq n
$.
\eit

\subsection{Main results}
A natural question (first posed in \cite{kont07-thesis}) is whether the conditions (P1)-(P3) completely characterize the
possible $(\etab_{ij})$ matrices, or if there are some other constraints that the $\eta$-mixing coefficients must satisfy. The main
technical result of this note is Theorem \ref{thm:mainhij}, which resolves this question in the affirmative. Thus, for any ``valid'' (i.e., 
satisfying (P1)-(P3)) $n \times n$ matrix 
$H=(h_{ij})$, there is 
a finite set $\X$ and a probability measure $\mu$ on $\X^n$ such that
$ \bar\h_{ij}(\mu) = h_{ij}$
for $1\leq i<j\leq n$.

More broadly, it is of interest to characterize the possible mixing rates that various processes may have. 
Chapter 26 of \cite{bradley07iii} deals with this question and gives several intricate constructions of random processes having prescribed
mixing rates, under various types of strong mixing. Following the work of Kesten and O'Brien \cite{kesten76}, it emerged that essentially arbitrary
mixing rates are possible for various mixing notions. Thus it is not surprising that the same holds true for $\eta$-mixing; this is an easy consequence
of our main result (Corollary \ref{cor:mixrate}).
% Theorem \ref{thm:mainhij}.

Along the way, we collect 
%other miscellaneous 
various other
observations regarding the $\eta$-mixing coefficients -- some of which are auxiliary in proving our main
results, and others 
may be of independent interest.
%given for future reference.
%which, while not particularly difficult or surprising,
%may prove useful if recorded in a readily referenced form.

\subsection{Notation}

We use the indicator variable
$\pred{\cdot}$ 
to assign 0-1 truth values 
to the predicate in 
$\set{\cdot}$. 

Random variables are capitalized ($X$), specified sequences
are written in lowercase ($x\in\X^n$), the shorthand
$\sseq{X}{i}{j}\defeq
\seq{X}{i}{j}
$ is used for all sequences, and
sequence concatenation is denoted multiplicatively:
$\sseq{x}{i}{j}\sseq{x}{j+1}{k}=\sseq{x}{i}{k}$.
Sums will range over the entire space of the summation variable;
thus
$\ds\sum_{\sseq{x}{i}{j}}f(\sseq{x}{i}{j})$ stands for
$$\ds\sum_{\sseq{x}{i}{j}\in
%\X^{j-i+1}
\sseq{\X}{i}{j}
}f(\sseq{x}{i}{j}),$$ 
where $\sseq{\X}{i}{j}$ is just $\X^{j-i+1}$, re-indexed for convenience.
\hide{
By convention, when $i>j$, we define
$$ \sum_{\sseq{x}{i}{j}}f(\sseq{x}{i}{j})
\equiv 
f(\eps)$$
where $\eps$ is the null sequence. }
For $y\in\sseq{\X}{1}{i}$ and $x\in\sseq{\X}{j}{n}$, we will write
$\mu(x\gn y)$ as a shorthand for
$\mu\set{\sseq{X}{j}{n}=x \gn \sseq{X}{1}{i}=y}$; no confusion should arise.

The {\em total variation} norm of a signed measure $\nu$ on $\X^n$ (i.e., vector 
%in 
$\nu\in\R^{\X^n}$) is defined by
\beq
\TV{\nu} = {\ts\oo2}\nrm{\nu}_1
= {\ts\oo2}\sum_{x\in\X^n} \abs{\nu(x)}
\eeq
(the factor of $1/2$ is not entirely standard).
Unless otherwise stated, $\X$ is a finite set.
Whenever we wish to be explicit about the dependence of $\h_{ij}$ and $\etab_{ij}$ on
a given measure $\mu$, we will write $\h_{ij}(\mu;y,w,w')$ and $\etab_{ij}(\mu)$, respectively.

\section{Constructions and proofs}
\hide{
-decreasing rows
-series-products $\oplus$
-growth rates O(n)
-parallel-products $\otimes$
-continuity
-consistent process, dependence on n
}

Let us begin with 
%a simple 
an easy
verification that (P3) holds for all $(\etab_{ij})$:
\belen
\label{lem:hijdecr}
Let $(\etab_{ij})_{1\leq i<j\leq n}$, be the 
%$n\times n$ 
$\h$-mixing matrix associated with a probability measure $\mu$ on $\X^n$.
Then, for all 
$1\leq i<j_1 < j_2\leq n$, we have
\beq
\bar\eta_{ij_2}\leq\bar\h_{ij_1}.
\eeq
\enlen
\bepf
Fix
%Pick
$1\leq i<j_1 < j_2\leq n$
and $y\in\sseq{\X}{1}{i-1},w,w'\in\sseq{\X}{i}{i}$. 
Then
\beq
\eta_{ij_2}(y,w,w') 
&=& 
{\ts\oo2}
\sum_{x\in\sseq{\X}{j_2}{n}}
\abs{
\mu(x \gn yw)
-    
\mu(x \gn yw')
}\\
&=&
{\ts\oo2}
\sum_{x\in\sseq{\X}{j_2}{n}}
\dsabs{
\sum_{u\in\sseq{\X}{j_1}{j_2-1}}
[
\mu(ux \gn yw)
-    
\mu(ux \gn yw')]
}\\
&\leq&
{\ts\oo2}
\sum_{x\in\sseq{\X}{j_2}{n}}
\sum_{u\in\sseq{\X}{j_1}{j_2-1}}
\abs{
\mu(ux \gn yw)
-    
\mu(ux \gn yw')
}\\
&=&
{\ts\oo2}
\sum_{z\in\sseq{\X}{j_1}{n}}
\abs{
\mu(z \gn yw)
-    
\mu(z \gn yw')
}\\
&=&
\eta_{ij_1}(y,w,w').
\eeq
\enpf

Next, we establish a simple continuity property of $\etab_{ij}$:
\belen
\label{lem:conth}
Suppose $\X$ is a finite set and
let $\calP_+^n(\X)$ be the set of all 
strictly positive
probability measures $\mu$ on $\X^n$ 
(i.e.,
%such that 
$\mu(x)>0$ for all $x\in\X^n$). 
Endow $\calP_+^n(\X)$ with the metric $\TV{\cdot}$. Then, for all $1\leq i<j\leq n$,
the functional
$\etab_{ij}: \calP_+^n(\X) \to \R$ is continuous with respect to $\TV{\cdot}$.
\enlen
\bepf
The continuity of 
$\h_{ij}(y,w,w'):\mu\mapsto\R$ 
for fixed $y\in\X^{i-1},w,w'\in\X$ follows immediately from
Lemma 5.4.1 of \cite{kont07-thesis}. The claim follows since 
continuity is preserved under finite maxima.
%the maximum of countably many continuous functions is continuous.
\enpf
\begin{rem}
%Note that 
Continuity breaks down on the boundary of $\calP_+^n(\X)$; see Section 5.4 of 
\cite{kont07-thesis} for an example.
\end{rem}

%The 
Our
construction of a measure with the desired mixing coefficients will proceed in stages, the final object being
composed of intermediate ones. The building blocks will be measures of a particular simple form.
For $1\leq k< n$, let $h\in\sseq{\R}{k+1}{n}$ be a vector of length $n-k$, satisfying 
$$0\leq h_{j+1} \leq h_j \leq 1 $$ 
for $k<j< n$; any such $h$ will be called a {\em valid $k$\th row}. We say that the measure $\mu$ on $\X^n$ is
{\em pure $k$\th row} (with respect to $h$) if its $\h$-mixing matrix 
$(\etab_{ij})_{1\leq i<j\leq n}$ satisfies
\beq
\etab_{ij} &= &
\pred{i=k}h_j.
\hide{
\left\{
\begin{array}{ll}
h_j
,& i=k \\
0
,& i\neq k.
\end{array}
\right.
}
\eeq

Our first technical result is the existence of arbitrary pure $k$\th row measures:
\belen
\label{lem:purek}
Fix $1\leq k<n$ and suppose
%Let 
$h\in\sseq{\R}{k+1}{n}$ 
is
%be
a valid $k$\th row vector. 
Then there exists a 
measure $\mu$ on
$\set{0,1}^n$
which is
pure $k$\th row with respect to $h$.
\enlen
\bepf
The proof will proceed by algorithmic construction. Let a valid 
$k$\th row vector
$h\in\sseq{\R}{k+1}{n}$ be given.
Initialize $\mu\supr{n+1}$ to be the uniform measure:
\beq
\mu\supr{n+1}(x) &=& 2^{-n},
\qquad x\in \set{0,1}^n.
\eeq
\hide{
For $1\leq i<j\leq n$, define $A_{ij}\subset \set{0,1}^n$ by
\beq
A_{ij} &=& \set{x\in\xoi^n : x_i = x_j}
\eeq
and let $B_{ij}=\xoi^n\setminus A_{ij}$.}
For $v\in[0,1]$, define the measure $\mu\supr{n,v}$ on $\xoi^n$ by
\beq
\mu\supr{n,v}(x) &=& \alpha_n(v)[v\pred{x_k=x_n}\mu\supr{n+1}(x) + (1-v)\pred{x_k\neq x_n}\mu\supr{n+1}(x)]
,
\eeq
where $\alpha_n(v)$ is the normalization constant ensuring that 
$\sum_x \mu\supr{n,v}(x)=1$,
and define $f_n:[0,1]\to[0,1]$ by
\beq
f_n(v) &=& \etab_{kn}(\mu\supr{n,v}).
\eeq
Lemma \ref{lem:conth} assures the continuity of $f_n$ and it is straightforward to verify that $f_n(0)=f_n(1)=1$ and $f_n(1/2)=0$.
Thus, there exists a $v^*\in[0,1]$ such that 
$
%\etab_{kn}(\mu\supr{n,v^*})
f_n(v^*)
=h_n$;
define the new measure $\mu\supr n$ by
\beqn
\label{eq:mun}
\mu\supr{n}(x) &=& \mu\supr{n,v^*}(x).
\eeqn
Similarly, for $v\in[0,1]$, define
\beq
\mu\supr{n-1,v}(x) &=& \alpha_{n-1}(v)[v\pred{x_k=x_{n-1}}\mu\supr{n}(x) + (1-v)\pred{x_k\neq x_{n-1}}\mu\supr{n}(x)],
\qquad x\in\xoi^n
\eeq
(where $\alpha_{n-1}(v)$ is again the appropriate normalization constant)
and define $f_{n-1}:[0,1]\to[0,1]$ by
\beq
f_{n-1}(v) &=& \etab_{k,n-1}(\mu\supr{n-1,v}).
\eeq
Again, it is easily seen that 
$f_{n-1}(0)=f_{n-1}(1)=1$ and $f_n(1/2)=h_n$, so by continuity there is a $v^*\in[0,1]$ for which
$
%\etab_{k,n-1}(\mu\supr{n-1,v^*})
f_{n-1}(v^*)
=h_{n-1}$;
so we may define the new measure
\beqn
\label{eq:mun1}
\mu\supr{n-1}(x) &=& \mu\supr{n-1,v^*}(x).
\eeqn
By construction, we have $\etab_{k,n-1}(\mu\supr{n-1})=h_{n-1}$; we claim that additionally,
\beqn
\label{eq:etabkn}
\etab_{k,n}(\mu\supr{n-1})=h_{n}
\eeqn
(in other words, the second modification 
%of the measure 
in (\ref{eq:mun1})
did not ``ruin'' the effects of the first modification
in (\ref{eq:mun})).
The 
%latter 
claim 
in (\ref{eq:etabkn})
holds because in fact for all $y\in\xoi^k$ and 
%all 
$x\in\xoi$, we have
\beqn
\label{eq:munn1}
\mu\supr{n}\set{X_n = x \gn \sseq{X}{1}{k}=y} 
&=&
\mu\supr{n-1}\set{X_n = x \gn \sseq{X}{1}{k}=y} ;
\eeqn
the latter fact is straightforward (though somewhat tedious) to verify.

We may now proceed by induction. Let $\mu\supr{t}$ be defined, for $k+1<t\leq n$. Define, for
$v\in[0,1]$,
\beq
\mu\supr{t-1,v}(x) &=& \alpha_{t-1}(v)[v\pred{x_k=x_{t-1}}\mu\supr{t}(x) + (1-v)\pred{x_k\neq x_{t-1}}\mu\supr{t}(x)],
\qquad x\in\xoi^n
\eeq
%(where $\alpha_{n-1}(v)$ is again the appropriate normalization constant)
and let $f_{t-1}:[0,1]\to[0,1]$ be
\beq
f_{t-1}(v) &=& \etab_{k,t-1}(\mu\supr{t-1,v}).
\eeq
Choose $v^*\in[0,1]$ so that $f_{t-1}(v^*)=h_{t-1}$ and define the new measure
\beq
\mu\supr{t-1} &\equiv \mu\supr{t-1,v^*}.
\eeq
Again, a straightforward calculation gives
\beqn
\label{eq:mutt}
\mu\supr{t}\set{\sseq{X}{t}{n} = x \gn \sseq{X}{1}{k}=y} 
&=&
\mu\supr{t-1}\set{\sseq{X}{t}{n} = x \gn \sseq{X}{1}{k}=y} 
\eeqn
for all
$y\in\xoi^k$ and all $x\in\xoi^{n-t+1}$, which ensures that
$$ 
\etab_{k,t-1}(\mu\supr{t-1}),
\etab_{k,t}(\mu\supr{t-1}),
\ldots,
\etab_{k,n}(\mu\supr{t-1})
$$
all have the right values.
The process terminates when we have constructed $\mu\supr{k+1}$; this is our desired pure $k$\th row measure with respect to $h$.
It remains to verify that $\etab_{ij}(\mu\supr{k+1})=0$ for $i\neq k$, but this is 
%straightforward.
%simple.
almost immediate.
\enpf
\begin{rem}
The ``backwards'' order of constructing the measures $\mu\supr{t}$ with $t=n,n-1,\ldots,k+1$ is essential. A construction
in the 
%``reverse''
``forward''
order fails precisely because (\ref{eq:mutt}) no longer holds. The reader is invited to verify that the marginals
of the constructed measure 
$\mu=\mu\supr{k+1}$
are identical, with $\mu\set{X_i=0}=\mu\set{X_i=1}=1/2$ for $1\leq i\leq n$.
\end{rem}

Next we turn to product measures. There are (at least) two natural ways to form products of probability measures; we shall refer to 
them as
{\em series} and {\em parallel}. Let $\calX,\calY$ be finite sets and $m,n \in \N$. 
If $\mu$ is a measure on $\calX^m$ and $\nu$ a measure on $\calX^n$, we define their series product, denoted by $\mu\serp\nu$,
to be the following measure on $\calX^{m+n}$:
\beqn
\label{eq:serp}
(\mu\serp\nu)(z) &=& \mu(x)\nu(y),
\qquad z=xy\in\calX^{m+n}, x\in\calX^m, y\in\calX^n.
\eeqn
If $\mu$ is a measure on $\calX^n$ and $\nu$ a measure on $\calY^n$, we define their parallel product, denoted by $\mu\parp\nu$,
to be the following measure on $(\calX\times\calY)^n$:
\beq
(\mu\parp\nu)(z) &=& \mu(x)\nu(y),
\qquad 
z=(x,y)\in (\calX\times\calY)^n.
\eeq

As our main construction will involve parallel products of measures, the following 
simple
result 
is 
useful.
%key (and may well be of independent interest):
\belen
\label{lem:parp}
Let $\mu$ and $\nu$ be probability measures on $\calX^n$ and $\calY^n$, respectively, and let
$\etab_{ij}(\mu)$, $\etab_{ij}(\nu)$ 
and $\etab_{ij}(\mu\parp\nu)$ 
be the corresponding $\eta$-mixing matrices.
Then we have
\beqn
\label{eq:parp}
\max\set{ \etab_{ij}(\mu),\etab_{ij}(\nu) } 
\;\leq\;
\etab_{ij}(\mu\parp\nu)
\;\leq\;
\etab_{ij}(\mu)+\etab_{ij}(\nu) 
%- \etab_{ij}(\mu)\etab_{ij}(\nu).
\eeqn
for all $1\leq i<j\leq n$.
\enlen
\bepf
Fix $i<j$. Throughout this proof, $x$ will denote sequences over $\calX$, $y$ sequences over $\calY$,
and $z=(x,y)$ over $\calX\times\calY$. Pick arbitrary 
$\sseq{z}{1}{i-1}=(\sseq{x}{1}{i-1},\sseq{y}{1}{i-1})$ and $z_i=(x_i,y_i)$, $z_i'=(x_i',y_i')$.
Then we expand
\beqn
\label{eq:parpcalc}
\h_{ij}(\mu\parp\nu;\sseq{z}{1}{i-1},z_i,z_i') &=&
\TV{ 
(\mu\parp\nu)(\cdot\gn \sseq{z}{1}{i-1} z_i) -
(\mu\parp\nu)(\cdot\gn \sseq{z}{1}{i-1} z_i')
} \\
\nonumber
&=& 
{\ts\oo2}
\sum_{\sseq{z}{j}{n}}
\abs{
(\mu\parp\nu)(\sseq{z}{j}{n}\gn \sseq{z}{1}{i-1} z_i) -
(\mu\parp\nu)(\sseq{z}{j}{n}\gn \sseq{z}{1}{i-1} z_i')
}
\\\nonumber
&=&
{\ts\oo2}
\sum_{\sseq{x}{j}{n}}
\sum_{\sseq{y}{j}{n}}
\abs{
\mu(\sseq{x}{j}{n}\gn \sseq{x}{1}{i-1} x_i) 
\nu(\sseq{y}{j}{n}\gn \sseq{y}{1}{i-1} y_i) 
-
\mu(\sseq{x}{j}{n}\gn \sseq{x}{1}{i-1} x_i') 
\nu(\sseq{y}{j}{n}\gn \sseq{y}{1}{i-1} y_i') 
}
\\\nonumber
&\geq&
{\ts\oo2}
\sum_{\sseq{x}{j}{n}}
\abs{
\sum_{\sseq{y}{j}{n}}
\sqprn{
\mu(\sseq{x}{j}{n}\gn \sseq{x}{1}{i-1} x_i) 
\nu(\sseq{y}{j}{n}\gn \sseq{y}{1}{i-1} y_i) 
-
\mu(\sseq{x}{j}{n}\gn \sseq{x}{1}{i-1} x_i') 
\nu(\sseq{y}{j}{n}\gn \sseq{y}{1}{i-1} y_i')
}
}
\\\nonumber
&=&
{\ts\oo2}
\sum_{\sseq{x}{j}{n}}
\abs{
\mu(\sseq{x}{j}{n}\gn \sseq{x}{1}{i-1} x_i) 
-
\mu(\sseq{x}{j}{n}\gn \sseq{x}{1}{i-1} x_i') 
}
\\\nonumber
&=&
\h_{ij}(\mu;\sseq{x}{1}{i-1},x_i,x_i').
\eeqn
Exchanging the roles of $x$ and $y$
%, this 
yields the lower bound in (\ref{eq:parp}).
To obtain the upper bound, we apply the $\TV{\cdot}$ tensorization property (see Lemma 2.2.5 in \cite{kont07-thesis}) to
(\ref{eq:parpcalc}):
\beq
&&
\TV{ 
(\mu\parp\nu)(\cdot\gn \sseq{z}{1}{i-1} z_i) -
(\mu\parp\nu)(\cdot\gn \sseq{z}{1}{i-1} z_i')
}
\;\leq\;
\\
&&
\TV{\mu(\cdot \gn \sseq{x}{1}{i-1} x_i) -\mu(\cdot \gn \sseq{x}{1}{i-1} x_i')}
+
\TV{\nu(\cdot \gn \sseq{y}{1}{i-1} y_i) -\nu(\cdot \gn \sseq{y}{1}{i-1} y_i')}
-
\\
&&
\TV{\mu(\cdot \gn \sseq{x}{1}{i-1} x_i) -\mu(\cdot \gn \sseq{x}{1}{i-1} x_i')}
\TV{\nu(\cdot \gn \sseq{y}{1}{i-1} y_i) -\nu(\cdot \gn \sseq{y}{1}{i-1} y_i')}
\eeq
which yields the desired bound.
\enpf

The interested reader may consult Lemma 3.2.1 of 
\cite{kont07-thesis} for some observations regarding the behavior of $\eta$-mixing coefficients under series products.

We are now ready to prove the main result of this note.
\bethn
\label{thm:mainhij}
Let $H=(h_{ij})$ be any $n \times n$ matrix satisfying (P1), (P2) and (P3).
Then there exists a finite set $\X$ and a probability measure $\mu$ on $\X^n$ such that
\beqn
\label{eq:mainhij}
\bar\h_{ij}(\mu) &=& h_{ij}
\eeqn
for $1\leq i<j\leq n$.
\enthn
\bepf
For $k=1,\ldots,n-1$, let $h\supr{k}\in\sseq{\R}{k+1}{n}$ be the vector $(h_{k,k+1},h_{k,k+2},\ldots,h_{k,n})$ -- i.e.,
the nonzero entries of the $k$\th row of $H$. Then Lemma \ref{lem:purek} provides a 
measure $\mu\supr{k}$ on $\xoi^n$
which is
pure $k$\th row with respect to $h\supr{k}$.
Let $\mu$ be the (parallel) product of these pure $k$\th row measures:
\beq
\mu &=& \mu\supr1 \parp\mu\supr2 \parp\ldots\mu\supr{n-1};
\eeq
note that $\mu$ is a measure on $\X^n$, where $\X=\xoi^{n-1}$. By definition of 
pure $k$\th row measures
and by Lemma \ref{lem:parp},
we have that (\ref{eq:mainhij}) holds.
\enpf

\begin{rem}
Our construction requires 
an exponential state space,
$\abs{\X}=2^{n-1}$. Are there analogous constructions using fewer states?
In Section 5.7 of \cite{kont07-thesis} we constructed a measure $\mu$ on $\xoi^n$ 
%having
satisfying (\ref{eq:mainhij}) for the special case where the rows of $H$ are constant:
$h_{i,i+1}=h_{i,i+2}=\ldots=h_{i,n}$; it seems unlikely that the general case is achievable with a
constant number of states.
\end{rem}

Up to this point, we have been discussing 
the $\eta$-mixing coefficients of
probability measures on finite sequences.
This notion extends quite naturally to random processes -- i.e., probability measures $\mu$ on $\X^\N$.
Let $\mu_n$ be the marginal distribution of $\sseq{X}{1}{n}$ and denote by
$\etab\supr{n}_{ij}$ the 
%$n\times n$ 
$\eta$-mixing matrix of $\mu_n$. It is straightforward to verify that in general,
%mixing coefficients
$\etab\supr{n}_{ij}$ depends on $n$ and that
\beq
\etab\supr{n}_{ij} &\leq& \etab\supr{n+1}_{ij}
\eeq
for $1\leq i<j\leq n$.
Let $\Delta_n(\mu)$ be the $n\times n$ matrix $\Delta$ corresponding to $\mu_n$, as defined in 
(\ref{eq:GDdef}). Recall that the $\ell_\infty$ operator norm of a nonnegative matrix is its maximal row sum.
Thus we can define the $\eta$-mixing rate of the process $\mu$ as the function $R_\mu:\N\to\R$:
\beq
R_\mu(n) &=& \nrm{\Delta_n(\mu)}_\infty.
\eeq
It's clear that
% \\
(i) $R_\mu$ is nondecreasing
%\\
and
%\\
(ii) $1\leq R_\mu(n)\leq n$;
%\\
any function satisfying these properties will be called a {\em valid} rate
function. 

\becon
\label{cor:mixrate}
Let $r:\N\to\N$ be a valid rate function. Then there is a set $\X=\xoi^\N$ and a measure $\mu$ on $\X^\N$ such that
\beqn
\label{eq:limsup}
\limsup_{n\to\infty} \frac{ R_\mu(n) }{ r(n) } &=& 1.
\eeqn
\encon
\bepf
We begin with the simple observation that if $r$ is a valid rate function then for all $k\geq1$ and all $0<\eps<1$,
there is an $n=n(k,\eps)>k$ and an $h=h(k,\eps)\in[0,1]$ such that
\beqn
\label{eq:nk}
1-\eps \;\leq\;
\frac{ h(k,\eps)(n-k) }{r(n) }
\;\leq\; 1.
\eeqn
Let $1>\eps_1>\eps_2>\ldots>0$ be a sequence decreasing to $0$. 
Pick a $k\geq1$ and let 
%$\eps=\eps_k$. 
%Let 
%$n_k=n(k,\eps_k)$ and $h=h(k,\eps_k)$, 
$n(k)=n(k,\eps_k)$ 
and 
$h(k)=h(k,\eps_k)$, 
as stipulated in (\ref{eq:nk}).
Define $h\supr{k}\in\sseq{\R}{k+1}{n}$ by
\beq
h\supr{k}_j &=& h(k),
\qquad k< j\leq n(k),
\eeq
and let $\mu\supr{k}$ be the measure on $\xoi^{n(k)}$
which is
pure $k$\th row
with respect to $h\supr k$, as constructed in Lemma \ref{lem:purek}.
Let $\beta$ be the symmetric Bernoulli measure on $\xoi$ (i.e., $\beta(0)=\beta(1)=1/2$) and define the 
measure $\hat\mu\supr{k}$ on $\xoi^\N$ by
\beq
\hat\mu\supr{k} &=& 
%\beta^{k-1}\serp 
\mu\supr{k} \serp \beta \serp \beta \serp \ldots
\eeq
where the operation $\serp$ is defined in (\ref{eq:serp}).
% and $\beta^{k-1}$ is the uniform measure on $\xoi^{k-1}$.
In this way, we have obtained a countable collection of measures
$\set{\hat\mu\supr{k} : k=1,2,\ldots}$ on $\xoi^\N$; note that 
by construction,
we have for each $k$
\beqn
\label{eq:muk}
1-\eps_k \;\leq\;
\frac{ 
\nrm{\Delta_{n(k)}( \hat\mu\supr{k} )}_\infty
%h(k,\eps_k)(n(k)-k) 
}{r(n(k)) }
\;\leq\; 1.
\eeqn
Now let $\mu$ be the measure on $(\xoi^\N)^\N$ obtained by taking the (parallel) product of all the 
$\hat\mu\supr{k}$'s:
\beq
\mu &=& \hat\mu\supr{1} \parp \hat\mu\supr{2} \parp \ldots 
\eeq
(the $\parp$ operator is defined in (\ref{eq:parp})).
It remains to verify that $\mu$ is a well-defined probability measure on 
$\X^\N$, $\X=\xoi^\N$ by applying the
Ionescu Tulcea theorem (\cite[Theorem 6.17]{Kallenberg02}), 
and that
(\ref{eq:muk}) continues to hold when 
$\hat\mu\supr{k}$ is replaced with $\mu$
-- the latter is straightforward.
\footnote{
To accommodate infinite state spaces, the $\max$ in (\ref{eq:hbar}) needs to be replaced with $\sup$.
} 
\enpf
%Existence of processes, consistency, behavior of $\h_{ij}$, Ionescu-Tulcea
\begin{rem}
%Note that
Our construction required an uncountable state space, $\X=\xoi^\N$. Are analogous constructions possible with smaller $\X$?
Is there a construction achieving (\ref{eq:limsup}) with $\lim$ in place of $\limsup$?
\end{rem}

\section*{Acknowledgments}
My thesis advisor John Lafferty encouraged me to explore the question of measures having prescribed mixing coefficients. I also thank 
Gideon Schechtman for hosting and guidance at the Weizmann Institute. A special thanks to Richard Bradley for the 
very
helpful correspondence.

\bibliographystyle{plain}
\bibliography{../mybib}
%\bibliography{mybib}

\end{document}